\documentclass{amsart}
\usepackage{amscd}
\usepackage{amsfonts}
\usepackage{amsmath}
\usepackage{amssymb}
\usepackage{latexsym}
\usepackage{amsthm,color}
\usepackage[all]{xy}
\usepackage{epsfig}
\usepackage{graphicx}
\usepackage{color}
\usepackage{tikz}
\usetikzlibrary{cd}
\usetikzlibrary{calc}
\newtheorem{proposition}{Proposition}

\newtheorem{theorem}{Theorem}
\newtheorem{corollary}{Corollary}

\theoremstyle{definition}
\newtheorem{definition}{Definition}
\theoremstyle{definition}
\newtheorem{example}{Example}
\theoremstyle{definition}
\newtheorem{remark}{Remark}

\newtheorem{problem}{Problem}

\newcommand{\Z}{\mathbb{Z}}
\newcommand{\R}{\mathbb{R}}

\begin{document}
\title[
Envelopes of straight line families in the plane 
]
{Envelopes of straight line families in the plane 
}
\author[Takashi Nishimura]{Takashi Nishimura
}
\address{
Research Institute of Environment and Information Sciences,
Yokohama National University,
Yokohama 240-8501, Japan}
\email{nishimura-takashi-yx@ynu.ac.jp}
\begin{abstract}
There is a widespread method to represent the envelope when a given 
hyperplane family creates an envelope.     
However, one sometimes encounters cases when the widespread method 
fails to represent the desired envelope precisely, and is confused.   
At the same time, one wants to find a correct method to draw the envelope 
precisely.       
\par 
In this article, focused on straight line families in the plane,  
an easy to understand explanation is given 
on the recently discovered correct method 
to represent the envelope precisely.   
Moreover, it is explained when and why the widespread method fails 
to represent the precise shape of envelope as well.     
\end{abstract}
\subjclass[2020]{57R45, 58C25} 
\keywords{Straight line family, Envelope, Gauss mapping, Frontal,
Creative, Creator.}


\date{}

\maketitle

\section{Introduction\label{section1}}
We start from an elementary example.   
\begin{example}\label{example1}
Let $f: \mathbb{R}\to \mathbb{R}^2$ be the mapping  
defined by $f(t)=\left(t, \sin t\right)$.     
The regular curve $f$ 
gives a parametrization of the non-singular curve 
\[
\mathcal{C}=\left\{(X, Y)\in \mathbb{R}^2\, |\, Y=\sin X\right\}.   
\]   
The affine tangent line $L_t$ to $\mathcal{C}$ 
at a point $\left(t, \sin t\right)$ 
may be defined by 
\[
\left(X-t, Y-\sin t\right)\cdot \left(-\cos t, 1\right)=0, 
\]
where the dot in the center stands for the standard scalar product of 
two vectors $\left(X-t, Y-\sin t\right)$ and 
$\left(-\cos t, 1\right)$ in the vector space $\R^2$.  
Since the straight line family $\left\{L_t\right\}_{t\in \mathbb{R}}$ 
is 
the affine tangent line family to $\mathcal{C}$, we believe  
that the sine curve $\mathcal{C}$ must be 
an envelope of $\left\{L_t\right\}_{t\in \mathbb{R}}$. 
Thus, by using the widespread method to represent  
the envelope of $\left\{L_t\right\}_{t\in \mathbb{R}}$ 
({\color{black}in the case of Example \ref{example1}, the widespread method 
means just to calculate the concrete form of the discriminant set 
$\mathcal{D}$ defined 
as follows.   For more details on} the widespread method, 
for instance refer to \cite{brucegiblin}), 
we try to confirm that $\mathcal{C}$ is actually an envelope of 
$\left\{L_t\right\}_{t\in \mathbb{R}}$.     
Set 
\[
F\left(X, Y, t\right)=\left(X-t, Y-\sin t\right)\cdot \left(-\cos t, 1\right)=
-\cos t X+Y+t\cos t - \sin t.   
\]
We have the following.   
\begin{eqnarray*}
\mathcal{D} & = & 
\left\{(X, Y)\in \mathbb{R}^2\: \left|\: \exists t\in \mathbb{R} \mbox{ s.t. }
F(X, Y, t)=\frac{\partial F}{\partial t}(X, Y, t)=0\right.\right\} \\ 
{ } & = & 
\left\{(X, Y)\in \mathbb{R}^2\: \left|\: \exists t\in \mathbb{R} \mbox{ s.t. }
-\cos t X+Y+t\cos t-\sin t=\sin t(X-t)=0\right.\right\} \\ 
{ } & = & \left\{(X, Y)\in \mathbb{R}^2\: 
\left|\: Y=X-2k\pi\; (k\in \Z) \mbox{ or }Y=-X+(2k+1)\pi\; (k\in \Z)
\mbox{ or }Y=\sin X \right.\right\} \\ 
{ } & \supsetneqq & \mathcal{C}.   
\end{eqnarray*}
Faced with the fact that $\mathcal{C}$ is a proper subset of $\mathcal{D}$, 
we are confused.     
In addition, we want to know a correct method 
to represent $\mathcal{C}$ precisely.   
Let us review the correct method given in \cite{nishimura} in the case of 
this example.     
%
There are three steps.   The first step is to normalize the defining equation 
$F=0$.   That is to say, replace 
%
the defining equation $F=0$ with a new one $G=0$  
having the form 
\[
G(X, Y, t)=
X\cos\theta(t)+Y\sin\theta(t)-a(t).    
\]
Then, we have (for example)  
\[
G(X, Y, t)=
X\cos\theta(t)+Y\sin\theta(t)-a(t) 
=X\frac{-\cos t}{\sqrt{\cos^2 t +1}}+Y\frac{1}{\sqrt{\cos^2 t +1}}
-\frac{-t\cos t+\sin t}{\sqrt{\cos^2 t+1}}{\color{black},}       
\]    
{\color{black}which implies that 
\[
\cos\theta(t)=\frac{-\cos t}{\sqrt{\cos^2 t +1}},\; 
\sin\theta(t)=\frac{1}{\sqrt{\cos^2 t +1}},\;  
a(t)=\frac{-t\cos t+\sin t}{\sqrt{\cos^2 t+1}}.   
\]
}
The second step (the most important step) 
is to find a $C^\infty$ function $b: \mathbb{R}\to \mathbb{R}$ 
satisfying 
\[
\frac{d a}{d t}(t)=b(t)\frac{d\theta}{d t}(t).   
\]
Elementary calculations show   
\[
\frac{d a}{d t}(t)=\frac{\sin t\left(t+\cos t\sin t\right)}
{\left(\cos^2 t+1\right)^{\frac{3}{2}}}, 
\quad \frac{d\theta}{dt}(t)= \frac{-\sin t}{\cos^2 t+1}.       
\]
Hence, we have $b(t)=\frac{-\left(t+\cos t\sin t\right)}{\sqrt{\cos^2 t+1}}$.     
The final step is just to substitute 
\[
a(t)=\frac{-t\cos t+\sin t}{\sqrt{\cos^2t+1}}, \; 
b(t)=\frac{-\left(t+\cos t \sin t\right)}{\sqrt{\cos^2 t+1}}, \; 
\left(\cos\theta(t), \sin\theta(t)\right) = 
\left(\frac{-\cos t}{\sqrt{\cos^2 t+1}}, \frac{1}{\sqrt{\cos^2 t+1}}\right) 
\]
into 
\[
a(t)\left(\cos\theta(t), \sin\theta(t)\right)+ 
b(t)\left(-\sin\theta(t), \cos\theta(t)\right).   
\]
Then, we have the desired parametrization of $\mathcal{C}$ as follows.   
\begin{eqnarray*}
{ } & { } & 
a(t)\left(\cos\theta(t), \sin\theta(t)\right)+ 
b(t)\left(-\sin\theta(t), \cos\theta(t)\right) \\ 
{ } & = & 
\frac{-t\cos t+\sin t}{\sqrt{\cos^2t+1}}
\left(\frac{-\cos t}{\sqrt{\cos^2 t+1}}, \frac{1}{\sqrt{\cos^2 t+1}}\right) 
+ 
\frac{-\left(t+\cos t \sin t\right)}{\sqrt{\cos^2 t+1}}
\left(\frac{-1}{\sqrt{\cos^2 t+1}}, \frac{-\cos t}{\sqrt{\cos^2 t+1}}\right) \\ 
{ } & = & 
\left(\frac{t\cos^2 t-\sin t\cos t+t+\cos t\sin t}{\cos^2 t+1}, 
\frac{-t\cos t+\sin t+t\cos t+\cos^2 t\sin t}{\cos^2 t+1}\right) \\ 
{ } & = & 
\left(\frac{t\left(\cos^2 t+1\right)}{\cos^2 t+1}, 
\frac{\sin t\left(1+\cos^2 t\right)}{\cos^2 t+1}\right) \\ 
{ } & = & 
\left(t, \sin t\right)  \\ 
{ } & = & f(t).    
\end{eqnarray*} 
\end{example}
\par 
\smallskip 
Let $N$, $a: N\to \mathbb{R}^{n+1}$ and $\nu: N\to S^n$ be 
an $n$-dimensional $C^\infty$ manifold without boundary, 
a $C^\infty$ function and a $C^\infty$ mapping 
respectively, where 
$S^n$ is the unit sphere in $\mathbb{R}^{n+1}$.   
For any $x\in N$, set 
$H_{\left(\nu(x), a(x)\right)}
=\left\{X\in \mathbb{R}^{n+1}\; |\; X\cdot \nu(x)=a(x)\right\}$.    
The problems on envelopes of hyperplanes
were classically studied (for instance see \cite{history}).    
Nevertheless, it is a surprizing fact that until very recently, 
the basic problems on envelopes were 
still wrapped in mystery (for instance, problems 
in Problem \ref{problem1} below seemed to be unsolved).  
{\color{black}Under the circumstances, in 2022, [7] gave a complete answer 
to each of} 
the following basic problems on envelopes created by  
hyperplane families 
$\mathcal{H}=\left\{H_{\left(\nu(x), a(x)\right)}\right\}_{x\in N}$.  
\begin{problem}\label{problem1}
\begin{enumerate}
\item[(1)] {\bf [Existence Problem]} 
Find a necessary and sufficient condition for 
a given hyperplane family to create an envelope.   
\item[(2)] {\bf [Uniqueness Problem]}  
Suppose that a given hyperplane family creates an envelope.   
Then, find a necessary and sfficient condition for the envelope to be unique.  
\item[(3)] {\bf [Representation Problem]}  
Suppose that a given hyperplane family creates an envelope.   
Then, find a representing formula of the envelope.   
\end{enumerate}
\end{problem} 
This paper is an easy to understand expository article on the 
solutions to Problem \ref{problem1} proved in \cite{nishimura}.          
In order to concentrate on explaining 
the core part of the solutions given in \cite{nishimura}, 
$n=1$ is assumed hereafter in this article.    
Namely, all answers to the following problems  
are explained in this article.    
\begin{problem}\label{problem2}
\begin{enumerate}
\item[(1)] {\bf [Existence Problem]} 
Find a necessary and sufficient condition for 
a given straight line family in the plane $\R^2$ to create an envelope.   
\item[(2)] {\bf [Uniqueness Problem]}  
Suppose that a given straight line family in the plane $\R^2$ 
creates an envelope.   
Then, find a necessary and sfficient condition for the envelope to be unique.  
\item[(3)] {\bf [Representation Problem]}  
Suppose that a given hyperplane family in the plane $\R^2$ 
creates an envelope.   
Then, find a representing formula of the envelope.   
\end{enumerate}
\end{problem} 
All answers (Theorem \ref{theorem1}, 
Theorem \ref{theorem2} and Theorem \ref{theorem3}) 
to Problem \ref{problem2}   
explained in this article are easily applicable to any concrete 
straight line family.     
For the proofs of Theorem \ref{theorem1}, Theorem \ref{theorem2} 
and Theorem \ref{theorem3}, see \cite{nishimura}.    
{\color{black}
Theorem \ref{theorem1} and Theorem \ref{theorem2}
have an advantage to clarify the reason why 
(1) and (2) of Problem \ref{problem2} 
have been neglected so far     
(see Corollary \ref{corollary1} 
and Corollary \ref{corollary2}). 
Theorem \ref{theorem3} 
has advantages to yield  
a natural generalization of Cahn-Hoffman vector formula given in 
\cite{hoffmancahn} and to play an important role to obtain 
a frontal version of the classical Legendre involution  
(see Corollary \ref{corollary3} and Corollary \ref{corollary4}). 
Moreover, Theorem \ref{theorem1} has more advantage to show that 
the existence of a singular point of the Gauss mapping characterizes 
the failure of 
widespread method (see Corollary \ref{corollary5}). 
To the best of author's knowledge, 
all of Theorem \ref{theorem1}, Theorem \ref{theorem2}, 
Theorem \ref{theorem3}, Corollary \ref{corollary1}, Corollary \ref{corollary2}, 
Corollary \ref{corollary3}, Corollary \ref{corollary4} 
and Corollary \ref{corollary5} are new results.    
}        
\par 
\bigskip 
This paper is organized as follows.   
In Section 2, we review several definitions concerning 
envelopes created by straight 
line families.    Section 3, Section 4 and Section 5 are devoted to 
explain Theorem 1 (Answer to Existence Problem of Problem \ref{problem2}), 
Theroem \ref{theorem2} 
(Answer to Uniqueness Problem of Problem \ref{problem2})  and 
Theorem \ref{theorem3} 
(Answer to Representation Problem of Problem \ref{problem2}) 
respectively.   
Finally, in Section 6, it is explained when and why the widespread method fails.    
\section{Preliminaries}\label{section2}
Any straight line $L$ in the plane $\mathbb{R}^2$ may be defined 
as follows, where $\theta, a$ are real numbers.  
\[
L=\{(X, Y)\in \mathbb{R}^2\; |\; X\cos\theta+Y\sin\theta=a\}.      
\]
Any straight line family $\mathcal{L}$ 
in the plane $\mathbb{R}^2$ may be defined 
as $\mathcal{L}=\left\{L_{\left(\theta(t), a(t)\right)}\right\}_{t\in \mathbb{R}}$, 
where $\theta, a: \mathbb{R}\to \mathbb{R}$ are 
$C^\infty$ functions and $L_{\left(\theta(t), a(t)\right)}$ is 
a straight line as follows.    
\[
L_{\left(\theta(t), a(t)\right)}=
\{(X, Y)\in \mathbb{R}^2\; |\; X\cos\theta(t)+Y\sin\theta(t)=a(t)\}.       
\]
\begin{definition}\label{Gauss}
{\rm 
Given a straight line family 
$\mathcal{L}=\left\{L_{\left(\theta(t), a(t)\right)}\right\}_{t\in \mathbb{R}}$, 
the mapping $\nu: \mathbb{R}\to S^1$ defined by 
$\nu(t)= \left(\cos\theta (t), \sin\theta (t)\right)$ 
is called the \textit{Gauss mapping} of 
$\mathcal{L}$. 
}
\end{definition}
\begin{definition}\label{envelope}
{\rm 
Let $\mathcal{L}=\left\{L_{\left(\theta(t), a(t)\right)}\right\}_{t\in \mathbb{R}}$ 
be a straight line family 
in the plane $\mathbb{R}^2$.   
A $C^\infty$ mapping $f: \mathbb{R}\to \mathbb{R}^2$ is called 
an \textit{envelope created by}  
$\mathcal{L}$ 
if the following two hold for any $t\in \mathbb{R}$.  
\begin{eqnarray*}
{\rm (a)} & \quad & \frac{d f}{d t}(t)\cdot \nu(t)  =  0,  \\    
{\rm (b)} & \quad & f(t)\in L_{\left(\theta(t), a(t)\right)}. \\
\end{eqnarray*} 
}
\end{definition}
\noindent 
Thus, an envelope is a $C^\infty$ mapping giving a solution of the 
first order linear differential equation (a) with one constraint condition 
(b).     
\begin{example}
Let $L_t$ and $f: \mathbb{R}\to \mathbb{R}^2$ be as in Example \ref{example1}.   
Thus, 
\[
L_t=\left\{(X, Y)\in \mathbb{R}^2\, 
|\, -\cos t X+Y+t\cos t-\sin t=0\right\}, \quad 
f(t)=\left(t, \sin t\right).    
\]
Then, it is easily seen that conditions (a), (b) in Definition \ref{envelope} 
are satisfied.    Hence, by definition, 
$f$ must be an envelope of the straight line family 
$\left\{L_t\right\}_{t\in\mathbb{R}}$.    
\end{example}   
\begin{definition}\label{frontal}
{\rm 
A mapping $f: \R\to \mathbb{R}^{2}$ is called 
a \textit{frontal curve} if there exists a mapping $\nu: \R\to S^1$ 
such that the following equality holds for any $t\in \R$.  
\[
\frac{d f}{d t}(t)\cdot \nu(t) =0.   
\]
} 
The mapping $\nu: \R\to S^1$ given above is called 
the \textit{Gauss mapping} of the frontal $f$.   
\end{definition}   
By definition, any envelope created by a straight line family is a frontal curve.   
Conversely, again by definition, 
any frontal curve $f: \R\to \mathbb{R}^{2}$ 
with Gauss mapping $\nu: \R\to S^1$ 
is an envelope of the straight line family 
$\mathcal{L}=\left\{L_{\left(\theta(t), a(t)\right)}\right\}_{t\in \mathbb{R}}$, 
where $\nu(t)=\left(\cos\theta(t), \sin\theta(t)\right)$ and  
$a(t)=f(t)\cdot \nu(t)$.     
Hence, these two notions are essentially the same 
although the notion of frontal curve has only 
recently been recognized and investigated.        
As an excellent survey article on frontal curve, 
\cite{ishikawa} is recommended to readers.    
\section{Answer to the existence problem on envelopes}\label{section3}
The following is the key notion in this article.   
\begin{definition}[\cite{nishimura}]\label{creative}
{\rm 
A straight line family 
$\mathcal{L}=\left\{L_{\left(\theta(t), a(t)\right)}\right\}_{t\in \mathbb{R}}$ 
in the plane $\mathbb{R}^2$ is said to be \textit{creative} if there exists 
a $C^\infty$ function $b: \R\to \R$ satisfying 
{\Large
\[
\underline{\frac{d a}{d t}(t)=b(t)\frac{d\theta}{d t}(t)}    
\leqno{(*)}
\]
}
for any $t\in \R$.   
The function $b:\R\to \R$ is called a \textit{creator}.    
}
\end{definition} 
\begin{theorem}[\cite{nishimura}]\label{theorem1}
A straight line family 
$\mathcal{L}=\left\{L_{\left(\theta(t), a(t)\right)}\right\}_{t\in \mathbb{R}}$ 
in the plane $\mathbb{R}^2$ creates an envelope if and only if it is creative.   
\end{theorem}
\begin{remark}
\begin{enumerate}
\item[(1)] In Example \ref{example1}, 
both $\cos\theta(t)$ and $\sin\theta(t)$ are 
even functions while $b(t)$ is an odd function. Hence,  
$b(t)$ cannot be 
obtained as the pullbuck 
(by $\nu$) 
of a $C^\infty$ section of the cotangent bundle 
$T^*S^1\to S^1$. 
\item[(2)] In Example \ref{example1}, we already calculated 
$\frac{d\theta}{dt}(t)=\frac{-\sin t}{\cos^2 t+1}$.    
This implies that $t=k\pi$ $(k\in \Z)$ are singular points of 
the Gauss mapping $\nu: \mathbb{R}\to S^1$ in Example 
\ref{example1}.   
Thus, as one can find in p.492 of \cite{wall}, 
for each $k\in \Z$ the dimension of the quotient vector space 
$V_k(\nu)/w\nu\left(V_k(1)\right)$ is greater than or equal to 1,  
where $V_k(\nu)$ is the vector space consisting of cotangent vector 
field germs $\xi: (\R, k\pi) \to T^*S^1$ along $\nu: (\R, k\pi)\to S^1$ 
and $w\nu\left(V_k(1)\right)$ is 
the vector space consisting of composition germs  
$\eta\circ \nu: (\R, k\pi)\to T^*S^1$ of the Gauss mapping germ 
$\nu: (\R, k\pi)\to S^1$ with $C^\infty$ section germs   
$\eta: (S^1, \nu(k\pi))\to T^*S^1$.      
Therefore, the above fact \lq\lq $b(t)$ cannot be obtained as the pullbuck 
of a $C^\infty$ section of the cotangent bundle 
$T^*S^1\to S^1$ by $\nu$\rq\rq\,  is not a surprizing fact.  
\item[(3)] By the above fact \lq\lq $b(t)$ cannot be obtained as the pullbuck 
of a $C^\infty$ section of the cotangent bundle $T^*S^1\to S^1$ 
by $\nu$\rq\rq\,,  it seems that Contact Geometry is useless for 
the proof of Theorem \ref{theorem1}.    Theorem \ref{theorem1} is proved by 
the anti-orthotomic technique developed in \cite{janeczkonishimura}.       
Moreover, not only Theorem \ref{theorem1} but also Theorem \ref{theorem2} in 
Section \ref{section4} and Theorem \ref{theorem3} in Section \ref{section5} 
can be obtained by the anti-orthotomic technique at once.    
\item[(4)] 
{\color{black}Just as the special case ($n=1$ case) of \cite{nishimura} 
Remark 1.1(a), it is reasonable to say that the height function $a: \R\to \R$ is 
\textit{differentiable with respect to} the Gauss mapping 
$\nu: \R\to S^1$ if there exists $b:\R\to \R$ satisfying 
$(*)$ in Definition \ref{creative}.        
Under this terminology, 
} 
the creator $b: \R\to \R$ may be 
called as the \textit{derived function} of the height function $a$ 
with respect to 
differentiation by the Gauss mapping $\nu$.    Hence, 
Theorem \ref{theorem1} may be replaced with the assertion 
\textit{\lq\lq A straight line family 
$\mathcal{L}=\left\{L_{\left(\theta(t), a(t)\right)}\right\}_{t\in \mathbb{R}}$ 
in the plane $\mathbb{R}^2$ creates an envelope 
if and only if the height function 
$a: \R\to \R$ is differentiable {\color{black}with respect to}  the Gauss mapping 
$\nu: \R\to S^1$\rq\rq}.    
\end{enumerate}      
\end{remark}
\begin{example}\label{example3}
{\rm 
Set $\theta(t)=a(t)= 0$.  
In this case, 
$\frac{d \theta}{d t}(t)=\frac{d a}{dt}(t)= 0$.  
Hence, an arbitrary $C^\infty$ function $b: \R\to \R$ is a creator, 
that is to say, it 
satisfies the equality 
$\frac{d a}{d t}(t)=b(t)\frac{d \theta}{d t}(t)\quad (\forall t\in \mathbb{R}).    
$
Therefore,  the straight line family 
$\{L_{\left(\theta(t), a(t)\right)}\}_{t\in \mathbb{R}}$ creates an envelope.   
}
\end{example}
\begin{example}\label{example4}
{\rm 
Set $\theta(t)= 0$,  $a(t)=t$.  
In this case, 
$\frac{d \theta}{d t}(t)= 0$ and 
$\frac{d a}{d t}(t)= 1$.   
Hence, there does not exist  $b: \mathbb{R}\to \mathbb{R}$ such that 
$\frac{d a}{d t}(t)=b(t)\frac{d \theta}{d t}(t)\quad (\forall t\in \mathbb{R}).    
$
Therefore,  the straight line family 
$\{L_{\left(\theta(t), a(t)\right)}\}_{t\in \mathbb{R}}$ is not creative.   
}
\end{example}
\begin{example}\label{example5}
{\rm 
Set $\theta(t)= t$,  $a(t)= 0$.  
In this case, 
$\frac{d \theta}{d t}(t)= 1$ and 
$\frac{d a}{d t}(t)= 0$.   
Hence, 
\[
0=\frac{d a}{d t}(t)=
\left(\frac{\frac{d a}{d t}(t)}{\frac{d\theta}{d t}(t)}\right)\frac{d \theta}{d t}(t) 
= \frac{0}{1}\times 1
\]
holds.
}
Thus, by Theorem \ref{theorem1}, 
the straight line family 
$\left\{L_{\left(\theta(t), a(t)\right)}\right\}_{t\in \mathbb{R}}$ creates 
an envelope.  
\end{example}
More generally, for a given straight line family 
$\mathcal{L}=\left\{L_{\left(\theta(t), a(t)\right)}\right\}_{t\in \mathbb{R}}$,  
suppose that the Gauss mapping 
$\nu: \R\to S^1$ of $\mathcal{L}$ is non-singular.  
Then, we have 
\[
\frac{d a}{dt}(t)=
\left(\frac{\frac{d a}{d t}(t)}{\frac{d\theta}{d t}(t)}\right)\frac{d \theta}{d t}(t) 
\]
for any $t\in \R$.   Thus, we have the following.   
\begin{corollary}\label{corollary1}
Let $\mathcal{L}=\left\{L_{\left(\theta(t), a(t)\right)}\right\}_{t\in \mathbb{R}}$ 
be a straight line family in the plane $\R^2$.    
Suppose that the Gauss mapping 
$\nu: \R\to S^1$ of $\mathcal{L}$ is non-singular.       
Then, the family $\mathcal{L}$ is always creative.    
Therefore, by Theorem \ref{theorem1}, it always creates an envelope.      
\end{corollary} 
\begin{example}\label{example6}
{\rm 
Set $\theta(t)= t^2$,  $a(t)= 0$.  
In this case, 
$\frac{d \theta}{d t}(t)=2t$ and 
$\frac{d a}{d t}(t)\equiv 0$.   
Thus, $t=0$ is a singular point of the Gauss mapping $\nu: \R\to S^1$.    
Nevertheless, we have 
\[
0=\frac{d a}{d t}(t)=0\times \frac{d\theta}{dt}(t).  
\]
for any $t\in \R$.     
}
Thus, by Theorem \ref{theorem1}, 
the straight line family 
$\left\{L_{\left(\theta(t), a(t)\right)}\right\}_{t\in \mathbb{R}}$ creates 
an envelope.  
\end{example}
\begin{example}[Evolute of the graph of $Y=\sin X$]\label{evolute}
{\rm 
Consider 
the affine normal line to the graph of 
$Y=\sin X$ at $\left(t, \sin t\right)$}.   
The defining equation of it may be  
\[
F\left(X, Y, t\right)=\left(X-t, Y-\sin t\right)\cdot 
\left(1, \cos t\right)=X+Y\cos t -t-\cos t\sin t=0.
\]
Thus, the normalized defining equation $G(X, Y, t)=0$ may be 
\[
G\left(X, Y, t\right)=
X\frac{1}{\sqrt{1+\cos^2 t}}+Y\frac{\cos t}{\sqrt{1+\cos^2 t}}
-\frac{t+\cos t\sin t}{\sqrt{1+\cos^2 t}}=0.
\]
Hence, we have 
$(\cos\theta(t), \sin\theta(t))=
\left(\frac{1}{\sqrt{1+\cos^2 t}}, \frac{\cos t}{\sqrt{1+\cos^2 t}}\right)$ and 
$a(t)=\frac{t+\cos t\sin t}{\sqrt{1+\cos^2t}}$.   
By calculations, we have   
\[ 
\frac{d a}{dt}(t)
=\frac{\cos t\left(3\cos t+\cos^3 t +t\sin t\right)}
{\left(1+\cos^2t\right)^{\frac{3}{2}}}, \quad    
\frac{d \theta}{dt}(t)=\frac{-\sin t}{1+\cos^2 t}.
\] 
For any $k\in \Z$, we have $\frac{da}{dt}(k\pi)\ne 0$ 
and $\frac{d\theta}{dt}(k\pi)=0$.    
Thus, there are no creator $b: \R\to \R$.  
Therefore, by Theorem \ref{theorem1}, 
the straight line family 
$\left\{L_{\left(\theta(t), a(t)\right)}\right\}_{t\in \R}$ 
does not create an envelope, namely, 
the evolute of the graph of $Y=\sin X$ does not exist.     
\end{example}
\section{Answer to the uniqueness problem on envelopes}\label{section4}
\begin{theorem}[\cite{nishimura}]\label{theorem2}
Let $\mathcal{L}=\left\{L_{\left(\theta(t), a(t)\right)}\right\}_{t\in \mathbb{R}}$ 
be a straight line family in the plane $\R^2$.    
Suppose that $\mathcal{L}$ is creative.   
Then, $\mathcal{L}$ creates a unique envelope if and only if 
the set consisting of regular points of the Gauss mapping 
$\nu: \R\to S^1$ of $\mathcal{L}$ is dense in $\R$.   
\end{theorem}
\begin{corollary}\label{corollary2}
Let $\mathcal{L}=\left\{L_{\left(\theta(t), a(t)\right)}\right\}_{t\in \mathbb{R}}$ 
be a straight line family in the plane $\R^2$.    
Suppose that the Gauss mapping 
$\nu: \R\to S^1$ of $\mathcal{L}$ is non-singular.       
Then, the family $\mathcal{L}$ creates a unique envelope.      
\end{corollary}
\begin{example}\label{example8}
As in Example \ref{example1}, set 
\[
\left(\cos\theta(t), \sin\theta(t)\right) = 
\left(\frac{-\cos t}{\sqrt{\cos^2 t+1}}, \frac{1}{\sqrt{\cos^2 t+1}}\right), \; 
a(t)=\frac{\sin t - t\cos t}{\sqrt{\cos^2 t+1}}.   
\] 
Then, in Example \ref{example1}, we confirmed 
$\frac{d\theta}{dt}(t)= \frac{-\sin t}{\cos^2 t+1}$.
Thus, the set of regular points of the Gauss mapping is a dense set  
in $\R$.   
Therefore, by Theorem \ref{theorem2}, the graph of $Y=\sin X$ must be the 
unique envelope.      
\end{example}
\begin{example}\label{example9}
{\rm 
As in Example \ref{example3}, set $\theta(t)=a(t)= 0$.  
In Example \ref{example3}, we confirmed that 
the straight line family 
$\{L_{\left(\theta(t), a(t)\right)}\}_{t\in \mathbb{R}}$ creates an envelope by 
Theorem \ref{theorem1}.    
Since $\frac{d \theta}{d t}(t)= 0$, the set of regular points of 
the Gauss mapping is empty.   Thus, by Theorem \ref{theorem2}, 
the created envelopes are not unique.   
}
\end{example}
\begin{example}\label{example10}
{\rm 
As in Example \ref{example6}, set $\theta(t)=t^2, \;a(t)=0$.  
We already confirmed that 
the straight line family 
$\{L_{\left(\theta(t), a(t)\right)}\}_{t\in \mathbb{R}}$ creates an envelope by 
Theorem \ref{theorem1}.    
Since $\frac{d \theta}{d t}(t)=2t$, the set of regular points of 
the Gauss mapping is dense in $\R$.   Thus, by Theorem \ref{theorem2}, 
the created envelope is unique.   
}
\end{example}
\section{Answer to the representation problem on envelopes}\label{section5}
\begin{figure}[h]
\begin{center}
\includegraphics[width=10cm]
{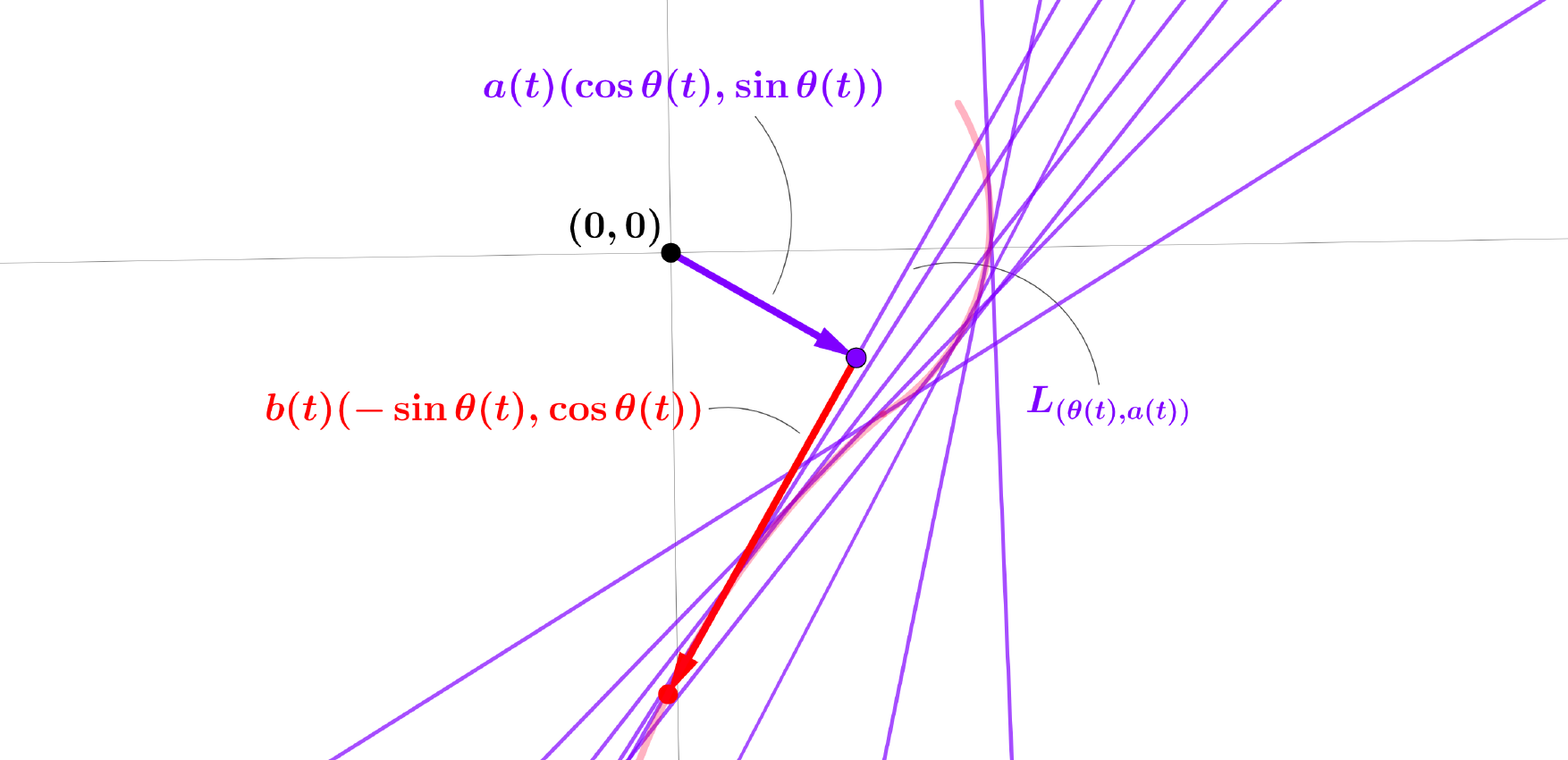}
\caption{Figure for Theorem \ref{theorem3}
}
\label{figure1}
\end{center}
\end{figure}

\begin{theorem}[\cite{nishimura}]\label{theorem3}
Let $\mathcal{L}=\left\{L_{\left(\theta(t), a(t)\right)}\right\}_{t\in \mathbb{R}}$ 
be a straight line family in the plane $\R^2$.    
Suppose that $\mathcal{L}$ is creative.   
Then, any envelope of $\mathcal{L}$ is parametrized by the mapping   
{\Large 
\[
\underline{
\R\ni t \mapsto a(t)\left(\cos\theta(t), \sin\theta(t)\right) + 
b(t)\left(-\sin\theta(t), \cos\theta(t)\right)\in \R^2
}, 
\]
}
where, $b: \R\to \R$ is the 
creator defined by $(*)$ of Definition \ref{creative}.     
\end{theorem} 
\noindent 
Theorem \ref{theorem3} is depicted in Figure \ref{figure1}.   
\begin{corollary}\label{corollary3}
Let $\mathcal{L}=\left\{L_{\left(\theta(t), a(t)\right)}\right\}_{t\in \mathbb{R}}$ 
be a straight line family in the plane $\R^2$.    
Suppose that the Gauss mapping 
$\nu: \R\to S^1$ of $\mathcal{L}$ is non-singular.       
Then, the unique envelope created by the family $\mathcal{L}$ 
is parametrized by the following mapping.  
\[
\R\ni t \mapsto a(t)\left(\cos\theta(t), \sin\theta(t)\right) + 
\left(\frac{\frac{d a}{dt}(t)}{\frac{d\theta}{dt}(t)}\right)
\left(-\sin\theta(t), \cos\theta(t)\right)\in \R^2.   
\]
\end{corollary}
There are two typical cases where the Gauss mapping is non-singular 
as follows.   
\begin{example}[Hedgehogs]\label{hedgehog}
{\rm 
Let $\mathcal{L}=
\{L_{\left(t, a(t)\right)}\}_{t\in \mathbb{R}}$ be a line family, where  
$a : \mathbb{R}\to \mathbb{R}$ is an arbitrary $C^\infty$ periodic 
function with period $2\pi$.    
Then, the Gauss mapping is non-singular.   
The {unique envelope} of $\mathcal{L}$ is called a \textit{hedgehog} 
which had been started to study by \cite{hedgehog}.    
In this case, since $\frac{d \theta}{d t}(t)=1$, the function $b: \R\to \R$ 
satisfying $(*)$ of Definition \ref{creative} is nothing but $\frac{d a}{d t}$.   
Thus, the hedgehog is parametrized as follows.     
\[
\R\ni t \mapsto a(t)\left(\cos t, \sin t\right) + 
\frac{d a}{d t}(t)\left(-\sin t, \cos t\right)\in \R^2.   
\] 
Therefore the celebrated {\color{black}C}ahn-Hoffman 
vector formula (\cite{hoffmancahn}) 
can be naturally obtained by our method.    
}
\end{example}
\begin{example}[Clairaut differential equations]
{\rm 
Consider a \textit{Clairaut differential equation}  
\[
Y=X\frac{dY}{dX}+g\left(\frac{dY}{dX}\right)
\leqno{(**)}
\]
where $g: \mathbb{R}\to \mathbb{R}$ is an arbitrary $C^\infty$ function.   
For any $t\in \mathbb{R}$, 
its {general solution $Y=tX+g(t)$} defines the straight line 
{$L_{\left(\theta(t), a(t)\right)}$},     
where 
$\nu(t)=\left(\cos\theta(t), \sin\theta(t)\right) 
=
\left(\frac{t}{\sqrt{t^2+1}},\frac{-1}{\sqrt{t^2+1}}\right)$ 
and $a(t)=\frac{-g(t)}{\sqrt{t^2+1}}$.    
It is easily seen that the Gauss mapping $\nu: \R\to S^1$ 
is non-singular.    
Thus, by Theorem \ref{theorem1} and Theorem \ref{theorem2}, 
there must exist the unique singular solution of 
the Clairaut differential equation $(**)$
as 
the unique envelope of the straight line family 
$\left\{L_{\left(\theta(t), a(t)\right)}\right\}_{t\in \R}$.   
Calculation shows that the unique creator $b: \R\to \R$ in this case has 
the following form.   
\[
b(t)=\frac{-\frac{dg}{dt}(t)\left(t^2+1\right)+tg(t)}{\sqrt{t^2+1}}.   
\]
Therefore, by Theorem \ref{theorem3}, the unique singular solution 
is parametrized as follows.   
\begin{eqnarray*}
{ } & { } & 
a(t)\left(\cos\theta(t), \sin\theta(t)\right) 
+ b(t)\left(-\sin\theta(t), \cos\theta(t)\right) \\ 
{ } & = & 
\frac{-g(t)}{\sqrt{t^2+1}}\left(\frac{t}{\sqrt{t^2+1}}, \frac{-1}{\sqrt{t^2+1}}\right)
+ 
\frac{-\frac{dg}{dt}(t)\left(t^2+1\right)+tg(t)}{\sqrt{t^2+1}}
\left(\frac{1}{\sqrt{t^2+1}}, \frac{t}{\sqrt{t^2+1}}\right) \\ 
{ } & = & 
\frac{1}{t^2+1}\left(-\frac{dg}{dt}(t)\left(t^2+1\right), 
\left(g(t)-t\frac{dg}{dt}(t)\right)\left(t^2+1\right)\right) \\ 
{ } & = & 
\left(-\frac{dg}{dt}(t),\; g(t)-t\frac{dg}{dt}(t)\right).    
\end{eqnarray*} 
For details on singularities of the unique singular solution of 
the Clairaut differential equation$(**)$, see for instance 
\cite{sajitakahashi}.   
} 
\end{example}
\begin{example}\label{example13}
{\rm 
As in Example \ref{example3} and Example \ref{example9}, 
set $\theta(t)=a(t)= 0$.  
We already confirmed that 
the straight line family 
$\{L_{\left(\theta(t), a(t)\right)}\}_{t\in \mathbb{R}}$ creates an envelope by 
Theorem \ref{theorem1} and the created envelopes are not unique 
by Theorem \ref{theorem2}.    
Let $b: \R\to \R$ be an arbitrary $C^\infty$ function.   
Since $b$ is a creator for envelopes of  
$\{L_{\left(\theta(t), a(t)\right)}\}_{t\in \mathbb{R}}$, by Theorem 
\ref{theorem3}, 
\[
a(t)\left(\cos\theta(t), \sin\theta(t)\right) + 
b(t)\left(-\sin\theta(t), \cos\theta(t)\right) 
 =  b(t)\left(-\sin 0, \cos 0\right)  
 =  \left(0, b(t)\right)
\]   
is an envelope of $\{L_{\left(\theta(t), a(t)\right)}\}_{t\in \mathbb{R}}$.    
}
\end{example}
\begin{example}\label{example14}
{\rm 
As in Example \ref{example6} and Example \ref{example10}, 
set $\theta(t)=t^2, \;a(t)=0$.  
We already confirmed that 
the straight line family 
$\{L_{\left(\theta(t), a(t)\right)}\}_{t\in \mathbb{R}}$ creates an envelope by 
Theorem \ref{theorem1} and the created envelope is unique  
by Theorem \ref{theorem2}.      
In this case, the unique creator is the constant function $0$.   
Therefore, by Theorem 
\ref{theorem3}, the unique envelope is parametrized as follows as desired.   
\[
 a(t)\left(\cos\theta(t), \sin\theta(t)\right) + 
b(t)\left(-\sin\theta(t), \cos\theta(t)\right) 
 =  \left(0, 0\right).    
\]    
}
\end{example} 
The following corollary, which may be regarded as the Legendre involution on  
frontals, can be obtained easily from Theorem \ref{theorem1} and 
Theorem \ref{theorem3}.    
\begin{corollary}\label{corollary4}
Suppose that an envelope 
\[
f(t)=a(t)\left(\cos\theta(t), \sin\theta(t)\right) 
+b(t)\left(-\sin\theta(t), \cos\theta(t)\right)
\]
of 
a straight line family 
$\{L_{\left(\theta(t), a(t)\right)}\}_{t\in \mathbb{R}}$ exists.   
Then, 
the mapping $g: \mathbb{R}\to \mathbb{R}^2$ defined by 
\[
g(t)=\left(b(t)\theta(t)-a(t)\right)
\left(\cos b(t), \sin b(t)\right) 
+\theta(t)\left(-\sin b(t), \cos b(t)\right)
\] 
is an envelope of the straight line family 
$\{L_{\left(b(t), \left(b(t)\theta(t)-a(t)\right)\right)}\}_{t\in \mathbb{R}}$.    
\end{corollary}  
\noindent 
\underline{Proof of Corollary4.}\quad By Theorem \ref{theorem1} and 
Theorem \ref{theorem3}, it follows that 
\[
\frac{d\left(b\theta\right)}{dt}(t)=b(t)\frac{d\theta}{dt}(t)
+\theta(t)\frac{d b}{dt}(t)=\frac{d a}{dt}(t)+\theta(t)\frac{d b}{dt}(t).
\]    
Thus, we have 
\[
\frac{d \left(b\theta-a\right)}{dt}(t)=\theta(t)\frac{d b}{dt}(t). 
\]   
Hence, again by Theorem \ref{theorem1} and Theorem \ref{theorem3}, 
the conclusion of Corollary \ref{corollary4} follows.    
\hfill 
$\Box$
\section{When and why the widespread method fails
}\label{section6}
Let $\mathcal{L}=
\left\{L_{\left(\theta(t), a(t)\right)}\right\}_{t\in \mathbb{R}}$ 
be a straight line family in the plane $\R^2$.    
Then, the normalized equation 
{
\[
G(X, Y, t)=X\cos\theta(t)+Y\sin\theta(t)-a(t)=0
\] 
}
is a {defining equation of $\mathcal{L}$}.     
{Suppose that $\mathcal{L}$ is creative}.   
Then, by Theorem \ref{theorem1}, we have 
\begin{eqnarray*}
\frac{\partial G}{\partial t}(X, Y, t) & = & 
\left(-X\sin\theta(t)+Y\cos\theta(t)\right)
\frac{d \theta}{d t}(t)-\frac{d a}{d t}(t) \\ 
{ } & = & 
\left(-X\sin\theta(t)+Y\cos\theta(t)
-b(t)\right)\frac{d \theta}{d t}(t).    
\end{eqnarray*}
Assume moreover that the Gauss mapping $\nu: \R\to S^1$ 
is non-singular.  
Then, {\color{black}by definition, 
it follows $\frac{d \theta}{d t}(t)\ne 0$ for any fixed $t\in \R$.   
Thus, for any fixed $t\in \R$, 
the Jacobian determinant of the linear mapping 
$\mathbb{R}^2\ni(X, Y)\mapsto 
\left(G(X, Y, t), \frac{\partial G}{\partial t}(X, Y, t)\right)\in\mathbb{R}^2$ 
is non-zero.   Hence, the solution of the system of linear equations $G(X, Y, t)=
\frac{\partial G}{\partial t}(X, Y, t)=0$ must be a unique $2$-dimensional vector.   
Therefore,} solving the system of equations $G(X, Y, t)=
\frac{\partial G}{\partial t}(X, Y, t)=0$ gives the precise shape of envelope of 
$\mathcal{L}$.    
Namely, the widespread method works well in this case 
to represent the envelope precisely.    
{\color{black}Hence we have the following.
\begin{proposition}\label{propositioin1}
Let  $\mathcal{L}=
\left\{L_{\left(\theta(t), a(t)\right)}\right\}_{t\in \mathbb{R}}$ 
be a straight line family in the plane $\R^2$.     
Suppose that the Gauss mapping of $\mathcal{L}$ is non-singular.   
Then, the unique envelope created by the straight line family can be 
directly obtained by the widespread method.       
\end{proposition}
}
\par 
Nextly, 
assume that the Gauss mapping $\nu$ is singular at  
$t=t_0$.     
{\color{black}
Then, by definition, we have 
\[
\frac{\partial G}{\partial t}(X, Y, t_0)=
\left(-X\sin\theta(t_0)+Y\cos\theta(t_0)
-b(t_0)\right)\frac{d \theta}{d t}(t_0)=0, 
\] 
which implies 
\[
\left\{(X, Y)\in \mathbb{R}^2\: \left|\: 
G(X, Y, t_0)=\frac{\partial G}{\partial t}(X, Y, t_0)=0\right.\right\}  
 =  
\left\{(X, Y)\in \mathbb{R}^2\: \left|\: 
G(X, Y, t_0)=0\right.\right\} 
 =  
L_{\left(\theta(t_0), a(t_0)\right)}.    
\]
}
Thus, the widespread method fails to get the precise shape of envelope 
in this case.   
{\color{black} 
Therefore, the following proposition holds.
}   
\begin{proposition}\label{proposition2}
Given a creative straight line family 
$\mathcal{L}=\left\{L_{\left(\theta(t), a(t)\right)}\right\}_{t\in \mathbb{R}}$ 
in the plane $\R^2$, suppose that 
the Gauss mapping of $\mathcal{L}$ is singular.     
Then, the widespread method fails to get an envelope precisely  
at and only at singular points of Gauss mapping.   
\end{proposition}   
{\color{black}
These two propositions yield }
\begin{corollary}\label{corollary5}
Given a creative straight line family 
$\mathcal{L}=\left\{L_{\left(\theta(t), a(t)\right)}\right\}_{t\in \mathbb{R}}$ 
in the plane $\R^2$, 
in order to get an envelope created by 
$\mathcal{L}$ precisely, 
the widespread method can be used when and only when 
the Gauss mapping is non-singular.   
\end{corollary}
\section*{{\color{black}Acknowledgements}}
{
{\color{black}The author wishes to express his sincere gratitude to the
referee for careful reading of this paper and making invaluable suggestions.
}
\par 
This work was partially supported 
by the Research Institute for Mathematical Sciences, 
a Joint Usage/Research Center located in Kyoto University. 
\par
The author was  
supported by JSPS KAKENHI (Grant No. 23K03109).   
}
%
%

\end{document}